\documentclass[12pt]{article}%
\usepackage[applemac]{inputenc}
\usepackage{graphicx}  
\usepackage{caption}   
\usepackage{amsmath,amssymb}
\usepackage{amsfonts}
\usepackage[applemac]{inputenc}
\usepackage{amsmath,amssymb}
\usepackage{color}
\usepackage{color, colortbl}
\usepackage{amsmath}
\usepackage{amssymb}
\usepackage{graphicx}
\usepackage{tikz}
\usepackage{geometry}
\usetikzlibrary{arrows}
\usepackage{amsfonts}
\usepackage{amsmath,amssymb}
\usepackage[applemac]{inputenc}
\usepackage{color}
\usepackage{amsmath}
\usepackage{amssymb}
\usepackage{graphicx}
\usepackage{tikz}
\newtheorem{theorem}{Theorem}[section]

\newtheorem{lemma}[theorem]{Lemma}
\newcommand{\E}{E}
\newtheorem{proposition}[theorem]{Proposition}
\newtheorem{corollary}[theorem]{Corollary}
\newtheorem{definition}[theorem]{Definition}

\newtheorem{example}[theorem]{Example}

\newtheorem{remark}[theorem]{Remark}

\usetikzlibrary{arrows,shapes,automata,petri}
\newenvironment{myenumerate}{

\begin{enumerate}}{\end{enumerate}}
\newcommand{\dproof}{\noindent {Proof.} \quad}
\newcommand{\fproof}{\hfill $\square$ \bigskip}

\numberwithin{equation}{section}

\definecolor{LightCyan}{rgb}{0.88,1,1}
\def\cF{{\mathcal {F}}}

\def\cB{{\mathcal{B}}}

\def\1B{\text{1\!\!I}}

\def\l{\langle}
\def\<{\langle}
\def\>{\rangle}
\def\P{\mathbb{P}}
\def\R{\mathbb{R}}
\def\l{\lambda}
\def\S{\mathcal{S}}
\def\E{\mathbb{E}}

\def\F{\mathcal{F}}
\def\H{\mathbb{H}}
\def\Y{\mathcal{Y}}
\def\G{\mathbb{G}}

\def\realio{\mathbb{R}_{0}}

\begin{document}

\title{The stochastic heat inclusion with fractional time derivative driven by  time-space Brownian white noise}

\author{Olfa Draouil$^{1}$, Rahma Yasmina Moulay Hachemi$^{2}$ \& Bernt \O ksendal$^{3}$
}
\date{27 November 2025}
\maketitle 
\footnotetext[1]{Department of Mathematics, University of Tunis El Manar, Tunisia\newline
Email: olfa.draouil@fst.utm.tn}

\footnotetext[2]{%
Department of Mathematics, University of Saida, Algeria.\newline
Email: yasmin.moulayhachemi@yahoo.com}

\footnotetext[3]{%
Department of Mathematics, University of Oslo, Norway. \\
Email: oksendal@math.uio.no.}


\begin{abstract}
We study a  time-fractional stochastic heat inclusion driven by additive time-space Brownian noise given by: 
\begin{equation}\nonumber
 	 \frac{\partial^{\alpha}}{\partial t^{\alpha}}Y(t,x)\in \l \Delta Y(t,x)+G(t,x,Y(t,x))+\sigma W(t,x);\; (t,x)\in (0,\infty)\times \mathbb{R}^{d}.
 	 	\end{equation}
  	Here $d\in\mathbb{N}=\{1,2,...\}$ is the space dimension, $x \in \R^d$ denotes a point in space, $t\geq 0$ denotes the time.
The fractional time derivative is interpreted as the Caputo derivative of order $\alpha \in (0,2).$ We show the following: \\
a) If a solution exists, then it is a fixed point of a specific set-valued map.\\
b) Conversely, any fixed point of this map is a solution of the heat inclusion.\\
 A solution $Y(t,x)$ is called \emph{mild} if $\E[Y^2(t,x)] < \infty$ for all $t,x$. We show that the solution is mild if\\
$\alpha=1$ \& $d=1,$ \ or \ $\alpha \geq 1$ \& $d\in \{1,2\}$.
On the other hand, if $\alpha < 1$ we show that the solution is not mild for any space dimension $d$.
  \end{abstract}

\section{Introduction}
The purpose of this paper is to  study fractional stochastic heat inclusions which can be used to model the transport of substances in porous media subject to noise and uncertainty. The applications show how our theoretical framework directly addresses real world environmental challenges, the combination of fractional derivatives, set valued maps and Brownian noise captures the phenomena. We illustrate this by an  example in Section 5 with numerical computations.

\begin{example}
Roughly, a \emph{differential inclusion} is a  \emph{set-valued differential equation}. For the readers who are not familiar with differential inclusions, we give a simple example:
\vskip 0.2cm
Consider the problem to find all differentiable functions $f: (0,\infty) \to (0,\infty)$ such that
\begin{align}
f'(x) \in (0, f(x)); \  \forall x>0; \  f(0)=1.
\end{align}
To solve this differential inclusion we first divide by $f(x)$ and  get
\begin{align*}
&\frac{f'(x)}{f(x)} \in (0,1), \  \text{ i.e.}\\
&\frac{d}{dx} \ln(f(x))\in (0,1); \  \forall x>0.
\end{align*}
Integrating this we get, using the Aumann definition of integration of set-valued functions and that $f(0)=1$,
\begin{align*}
&\ln(f(x)) \in \int_0^x (0,1) dt = \Big\{\int_0^x g(t)dt; \ 0<g(\cdot)<1, g \text{ integrable} \Big\}.
\end{align*}
Therefore the general solution of this inclusion is
\begin{align*}
f(x) \in \Big\{ \exp(\int_0^x g(t)dt)\Big\}_{0<g(\cdot) < 1}.
\end{align*}
\end{example}

\begin{example}
Consider the problem to find all differentiable functions $f: (0,\infty) \to (0,\infty)$ such that
\begin{align}
f'(x) \in (0, \frac{1}{f(x)}); \  \forall x>0; \  f(0)=1.
\end{align}
To solve this differential inclusion we first multiply by $f(x)$ and  get
\begin{align*}
&f(x)f'(x) \in (0,1), \  \text{ i.e.}\\
&\frac{d}{dx}(\tfrac{1}{2} f^2(x))\in (0,1); \  \forall x>0.
\end{align*}
Integrating this we get, using the Aumann definition of integration of set-valued functions and that $f(0)=1$,
\begin{align*}
&\tfrac{1}{2}f^2(x) \in \frac{1}{2}+ \int_0^x (0,1) dt =\tfrac{1}{2}+ \Big\{\int_0^x g(t)dt\Big\}_{0<g(\cdot)<1, g \text{ integrable}}=\Big\{\tfrac{1}{2}+ \int_0^x g(t)dt\Big\}_{0<g(\cdot)<1, g \text{ integrable}}.
\end{align*}
Therefore the general solution of this inclusion is
\begin{align*}
f(x) \in \Big\{ (1+\int_0^x 2 g(t)dt)^{\tfrac{1}{2}}\Big\}_{0<g(\cdot) < 1, g \text{ integrable}}.
\end{align*}
\end{example}    

\begin{example}
 Consider the problem to find all differentiable functions $f: (0,\infty) \to (0,\infty)$ such that
\begin{align}
f'(x) \in (0, f^2 (x)); \  \forall x>0; \  f(0)=1.
\end{align}
To solve this differential inclusion we first divide by $f^2(x)$ and  get
\begin{align*}
&\frac{f'(x)}{f^2(x)} \in (0,1), \  \text{ i.e.}\\
&\frac{d}{dx}(-\frac{1}{f(x)})\in (0,1); \  \forall x>0.
\end{align*}
Integrating this we get, using the Aumann definition of integration of set-valued functions and that $f(0)=1$,
\begin{align*}
&-\frac{1}{f(x)} \in -1 + \int_0^x (0,1) dt =-1+ \Big\{\int_0^x g(t)dt\Big\}_{g \in \G},
\end{align*}
where $\G$ is the set of all (measurable) integrable functions $g:(0,\infty) \to (0,1)$.\\
This is equivalent to saying that there exists a function $g \in \G$ such that
\begin{align*}
    -\frac{1}{f(x)}=-1 + \int_0^x g(t) dt\ \text{ i.e. }\\
    f(x)= (1 - \int_0^x g(t)dt)^{-1} \text{ for all } x>0.
    \end{align*}
We conclude that the general solution of this inclusion is
\begin{align*}
f(x) \in \Big\{ (1-\int_0^x  g(t)dt )^{-1}\Big\}_{g \in \G_0},
\end{align*} 
where $\G_0=\Big\{g:[0,\infty) \to (0,1); \text{ integrable and }  \int_0^x g(t)dt < 1 \ \forall x \Big\}.$

\begin{remark}
    It is interesting to compare with the corresponding differential \emph{equation}:
    \begin{align}
f'(x)= f^2(x); \ f(0)=1.
    \end{align}
\end{remark}
The unique  solution is
\begin{align}
    f(x)= \frac{1}{1-x},
\end{align}
which explodes for $x=1$. Thus in the equation case there is no global solution.
\end{example}
\subsection{Mathematical details of the model}
 We now explain our model more precisely:\\
 Let
 $$B(t,x)=B(t,x,\omega); t\geq 0, x=(x_1, x_2, ... ,x_d) \in \R^d, \omega \in \Omega$$
 be the time-space Brownian sheet on a filtered probability space $(\Omega,\F=\{\F_{t,x}\}_{(t,x)\in \R_+\times \R^d},\P)$, where $\P$ is the law of $B$.\\
 The time-space white noise $W(t,x)$ is defined by (in the sense of distributions)
\begin{equation}
  	W(t,x)=W(t,x,\omega)=\frac{\partial}{\partial t}\frac{\partial^{d}B(t,x)}{\partial x_{1}...\partial x_{d}}. 
  \end{equation} 

We define $\F_{t,x}$ to be the filtration generated by the Brownian sheet values $$\{B(s,a); s\leq t, a \leq x\}.$$
Let $\S'= \S'(\R^{1+d})$ denote the space of tempered distributions. Let $Y(t,x)=Y(t,x,\omega); (t,x) \in \R_{+} \times \R^d, \omega \in \Omega$ be an $\S'$-valued process and let
\begin{align}
G(t,x,Y(t,x)): \mathbb{R}_+\times\mathbb{
R}^d\times\S' \mapsto \mathcal{B}(\S') 
\end{align}
be a set-valued map, where $\mathcal{B}(\S')$ is the set of all Borel subsets of $\S'$.
\begin{definition} 
We call a function $g$ a \emph{selector} of the set $G$ if $g(t,x) \in G(t,x,Y(t,x))$ for all $(t,x)$. We let $\G$ denote the set of all (measurable) selectors of $G$.
\end{definition}
\begin{remark}
 Note that even though the set $G$ may depend on $Y$ all the selectors of $G$ may still be bounded. 
For example, this is the case if we fix $\varphi \in \S $ and define $G$ to be the interval $(0, \min(1,\langle Y(t,x),\varphi\rangle^2)$.
\end{remark}
 \begin{remark}\label{rem1.7}In this paper, we work with distributions, and in order to define the integral of distributions, we introduce the Gelfand-Pettis integral. This means that if $h:\R^n \mapsto \S'$ we are defining the integral $\int_{\R^n}h(x)dx $ as the element in $\S'$ whose action on a stochastic test function $\varphi \in \S$ (the Schwartz space of rapidly decreasing smooth functions) is given by
    \begin{align}
        \langle \int_{\R^n}h(x)dx, \varphi \rangle = \int_{\R^n} \langle h(x),\varphi\rangle dx.
    \end{align}
    This defines $\int_{\R^n}h(x)dx$ as an element of $\S'$ if $\langle h(x),\varphi\rangle$ is Lebesgue integrable in $\R^n$, which will always be the case in our application to $W(t,x)$.\end{remark}

\begin{definition}
 In the following we let $\mathcal{Y}$ denote the set of all Gelfand-Pettis integrable $\S'$ -valued processes $Y(t,x)=Y(t,x,\omega);  (t,x,\omega)\in \R_+ \times \R^d \times \Omega$.

\end{definition}
We study the following \emph{fractional stochastic heat inclusion} in the unknown $Y \in \mathcal{Y}$:
\begin{equation}\label{heat1}
 	 \frac{\partial^{\alpha}}{\partial t^{\alpha}}Y(t,x)\in \l \Delta Y(t,x)+G(t,x,Y(t,x))+\sigma W(t,x);\; (t,x)\in (0,\infty)\times \mathbb{R}^{d}.
 	 	\end{equation}
  	Here $d\in\mathbb{N}=\{1,2,...\}$ is the space dimension, $x \in \R^d$ denotes a point in space, $t\geq 0$ denotes the time,  $\frac{\partial^{\alpha}}{\partial t^{\alpha}}$ is the Caputo derivative of order $\alpha \in (0,2)$, and $\l>0$, $\sigma \in \mathbb{R}$ are given constants,
   
  	\begin{equation}
  		\Delta Y =\sum_{j=1}^{d}\frac{\partial^{2}Y}{\partial x_{j}^{2}}(t,x)
  	\end{equation}
  is the Laplacian operator with respect to $x$.

The initial and boundary conditions are
\begin{align}
    Y(0,x)&=\delta(x)\text{ (the Dirac measure at  } x), \label{1.4}\\
    \lim_{|x| \rightarrow \infty}Y(t,x)&=0.\label{1.5}
\end{align}
\begin{remark}
 Since  $W (t, x)=W(t,x,\omega)$ is the derivative (in the weak sense) of the continuous Brownian sheet $B(t,x,\omega)$ it can be represented as an element of $\S'(\R^{1+d})$ for each $\omega$.
    
\end{remark}

\begin{definition}
    We say that $Y \in \Y$ is a \emph{solution} of the fractional stochastic inclusion \eqref{heat1} if there is a selector $g \in \G$ such that
    \begin{equation}\label{heat2}
 	 \frac{\partial^{\alpha}}{\partial t^{\alpha}}Y(t,x) = \l \Delta Y(t,x)+g(t,x)+\sigma W(t,x);\; (t,x)\in (0,\infty)\times \mathbb{R}^{d}.
 	 	\end{equation}
        in the sense of distribution.
\end{definition}

\begin{remark}\textbf{(Explanation of the model)}
    
In the classical case, i. e. when $\alpha=1$, this equation models the normal diffusion of heat in a random  or noisy medium, the noise being represented by the time-space white noise $W(t,x)$. 

 - When $\alpha >1$ the inclusion can be used to model \emph{superdiffusion or enhanced diffusion}, where the particles spread faster than in regular diffusion. For example, this may occur in some biological systems.
 
 - When $\alpha <1$ the inclusion models \emph{subdiffusion}, in which the times of travel of the particles are longer than in the standard case. Such a situation may occur in transport systems.

 In all cases the drift term $G(t,x,Y(t,x))$ represents a sink/source. It is set-valued because it models an uncertainty in the possible values it can have. 
\end{remark}
 We consider the inclusion \eqref{heat1} in the space $\Y$, and in Theorem \ref{heat} we show that a solution $Y\in \Y$ (if it exists) must be a fixed point of a corresponding set-valued map and, conversely, any fixed point of this map is also a solution of the fractional heat inclusion. 

 Indeed, the solution \(Y(t, x)\) that we obtained is, in general, distribution-valued. However, in some cases, the solution can be represented as an element of \(L^{2}(P)\). A natural question is therefore to determine under which conditions the solution \(Y(t, x)\) belongs to \(L^{2}(P)\).
\\
We define a solution $Y(t,x)$ to be \emph{mild} if $Y(t,x) \in L^2(\P)$ for all $t,x$. 
It is well-known that in the classical heat \emph{equation} case with $\alpha = 1$, the solution is mild if and only if the space dimension $d=1$. See e.g. \cite{Hu}. In Theorem \ref{th8.1} we show that if $\alpha \in (1,2)$ the solution is mild if $d=1$ or $d=2$, 
and we show that if $\alpha < 1$ then the solution is not mild for any space dimension $d$.\\

Here is an outline of the rest of our paper:
\begin{itemize}
\item 
In Section 2 we give some mathematical preliminaries. 
\item 
Section 3 contains our main result.
\item 
In Section 4 we study the question: When is the solution mild?
\item 
Finally, in Section 5 we present an example, where we include numerical calculations and figures.
 \end{itemize} 
 
\section{Mathematical preliminaries}
In this section we recall briefly some of the main concepts and background results we will need.
\subsection{The space of tempered distributions}
 We recall some of the basic properties of the Schwartz space $\mathcal{S}$ of rapidly decreasing smooth functions and its dual, the space $\mathcal{S}'$ of tempered distributions.

Let $n$ be a given natural number. Let $\mathcal{S}=\mathcal{S}(\mathbb{R}^n)$\label{simb-028} be the
space of rapidly decreasing smooth real
functions $f$
on $\mathbb{R}^n$
equipped with the family of seminorms:\label{simb-029} 
\begin{equation*}
\Vert f \Vert_{k,\gamma} := \sup_{x \in \mathbb{R}^n}\big\{ (1+|x|^k) \vert
\partial^ \gamma f(x)\vert \big\}< \infty,
\end{equation*}
where $k = 0,1,...$, $\gamma=(\gamma_1,...,\gamma_n)$ is a multi-index with $%
\gamma_j= 0,1,...$ $(j=1,...,n)$ and\label{simb-030} 
\begin{equation*}
\partial^\gamma f (x):= \frac{\partial^{|\gamma|}}{\partial
x_1^{\gamma_1}\cdots \partial x_n^{\gamma_n}}f(x)
\end{equation*}
for $|\gamma|=\gamma_1+ ... +\gamma_n$.

Then
$\mathcal{S}=\mathcal{S}(\mathbb{R}^n)$ is a
Fr\'echet space.

Let $\mathcal{S}^{\prime }=\mathcal{S}^{\prime }(\mathbb{R}^{n})$\label%
{simb-031} be its dual, called the space of \emph{tempered distributions}. 
\index{tempered distributions} Let $\mathcal{B}$ denote the family of all
Borel subsets of $\mathcal{S}^{\prime }(\mathbb{R}^{n})$ equipped with the
weak* topology. If $\Phi \in \mathcal{S}^{\prime }$ and $f \in \mathcal{%
S}$ we let \label{simb-033} 
\begin{equation}
\Phi (f) \text{ or } \langle \Phi ,f \rangle  \label{3.1}
\end{equation}%
denote the action of $\Phi$ on $f$.

\begin{example}
\begin{itemize}
\item
{(Evaluations)}
For $y \in \R$ define the function $\delta_y$ on $\S(\R)$ by $\delta_y(\phi)=\phi(y)$. Then $\delta_y$ is a tempered distribution.\vskip 0.2cm
\item
{(Derivatives)} Consider the function $D$, defined for $\phi \in \S(\mathbb{R})$ by $D[\phi]=\phi^{\prime}(y)$. Then  $D$ is a tempered distribution. \vskip 0.2cm
\item
{(Distributional derivative)}\\
 Let $T$ be a tempered distribution, i.e. $T \in \S^{'}(\mathbb{R}) $. We define the distributional derivative $T^{'}$ of $T$ by
 $$ T^{'}[\phi]=-T[\phi^{'}]; \quad \phi \in \S.$$
 Then $T^{'}$ is again a tempered distribution.
 \end{itemize}
 \end{example}

 \subsection{The Mittag-Leffler functions}

 The following functions are useful when dealing with fractional derivatives:
 \begin{definition}(The two-parameter Mittag-Leffler function)
 The Mittag-Leffler function of two parameters $\alpha,\; \beta$ is denoted by $E_{\alpha,\beta}(z)$ and defined by:
\begin{equation}
    E_{\alpha,\beta}(z)=\sum_{k=0}^{\infty}\frac{z^{k}}{\Gamma(\alpha k+\beta)}
\end{equation}
where $z,\; \alpha,\; \beta\in \mathbb{C},\; Re(\alpha)>0\; and\; Re(\beta)>0,$ and $\Gamma$ is the Gamma function.
\end{definition}

 \begin{definition}(The one-parameter Mittag-Leffler function)
The Mittag-Leffler function of one parameter $\alpha$ is denoted by $E_{\alpha}(z)$ and defined as;
\begin{equation}
    E_{\alpha}(z)=\sum_{k=0}^{\infty}\frac{z^{k}}{\Gamma(\alpha k+1)}
\end{equation}
where $z,\; \alpha\in \mathbb{C},\; Re(\alpha)>0.$
\end{definition}

\subsection{The Caputo fractional derivative}
  We summarize the definitions and some properties of the Caputo derivatives:
  \begin{definition}
  The Caputo fractional derivative of order $\alpha > 0$ of a function $f$ such that $f(x)=0$ when $x<0$ is denoted by  $D^{\alpha} f (x)$ or $\frac{d^{\alpha}}{dx^{\alpha}} f(x)$ and defined  by 
  \begin{align}\label{caputo1}
  D^{\alpha}f(x):& =
  \begin{cases}
  \frac{1}{\Gamma(n-\alpha)}\int_0^x \frac{f^{(n)}(u)du}{(x-u)^{\alpha +1 -n}}; \quad n-1 < \alpha < n\\
  \frac{d^n}{dx^n}f(x); \quad \alpha =n.
  \end{cases}
  \end{align}
  Here $n$ is an smallest integer greater than or equal to $\alpha$.\\
  
  \noindent If $f$ is not smooth these derivatives are interpreted in the sense of distributions, as explained in Section 2.1.
  \end{definition}
 
 \begin{example}
  If $f(x)=x$ and $\alpha \in (0,1)$ then
 
  \begin{align}
 D^{\alpha}f(x)=\frac{ x^{1-\alpha} }{(1-\alpha)\Gamma(1-\alpha)}.
 \end{align} 
In particular, choosing $\alpha=\tfrac{1}{2}$ we get
 \begin{align}
D^{\tfrac{1}{2}}f(x)=\frac{2 \sqrt{x}}{\sqrt{\pi}}.
\end{align} 
\end{example}

\subsection{Laplace transform of Caputo derivatives} 

Recall that the Laplace transform $L$ is defined by 
 
\begin{equation}
	Lf(s)=\int_{0}^{\infty}e^{-st}f(t)dt=:\widetilde{f}(s)
	\end{equation}
	for all functions $f$ such that the integral converges.\\
	Some useful properties of the Laplace transform and its relation to Mittag-Leffler functions and fractional derivatives, are the following:

  \begin{align}
    &L[ \frac{\partial ^{\alpha}}{\partial t^{\alpha}}f(t)](s)=s^{\alpha}(L f)(s)-s^{\alpha-1}f(0) \label{L1}\\       
     &L[E_{\alpha}(bt^{\alpha})](s)  = \frac{s^{\alpha -1}}{s^{\alpha}-b}\label{L2}\\
     &L[t^{\alpha-1}E_{\alpha,\alpha}(-b t^{\alpha})](s)=\frac{1}{s^{\alpha}+b}.\label{L3}
     \end{align}
     The convolution $f\ast g$ of two functions $f,g: [0,\infty) \mapsto \mathbb{R}$  is defined by
\begin{align}
(f \ast g)(t)=\int_0^t f(t-r)g(r) dr; \quad t \geq 0.
\end{align}

The convolution rule for Laplace transform states that $$L\left( \int_{0}^{t}f(t-r)g(r)dr\right) (s)=Lf(s)Lg(s),$$ 
or 
\begin{equation}\label{12}
	\int_{0}^{t}f(t-w)g(w)dw=L^{-1}\left( Lf(s)Lg(s)\right) (t).
\end{equation}

\subsection{Integrals of set-valued functions}
  To get a precise interpretation of our fractional differential inclusion \eqref{heat} we need to define what we mean by the integral of a set-valued function.  We will adopt the definition of Aumann \cite{Au}, as follows:
\begin{definition}
  Let $G: \R \mapsto \mathcal{B}$ be a set valued function; $-\infty \leq a < b \leq \infty.$
  Then we define the integral of $G(t)$ by
\begin{align}
  \int_a^b G(t)dt = \Big\{ \int_a^b g(t)dt\Big \}_{ g \in \mathbb{G}}
  \end{align}
  where $\mathbb{G}$ is the set of all integrable functions $g:[a,b] \mapsto \R$ such that $g(t) \in G(t)$ for all $t \in [a,b]$.
  \end{definition} 

\subsection{The Laplace transform and Fourier transform of a set-valued function}
Let $G$ be as above.
Then we define the Laplace transform of $G$ as follows:
\begin{definition}
\begin{align}
L(G)(s):&=\widetilde{G}(s):=\int_0^{\infty} e^{-st}G(t)dt \nonumber\\
&= \{ \int_0^{\infty} e^{-st} g(t)dt; g \in \mathbb{G} \}=\{ \widetilde{g}(s); g \in \mathbb{G} \};\ s\geq 0.
\end{align}
\end{definition}
The Fourier transform of $G$ is defined by
\begin{definition}
\begin{align}
F(G)(y):&=\widehat{G}(y):=\int_{\R^d} e^{-ixy}G(x)dx \nonumber\\
&= \{ \int_{\R^d} e^{-ixy} g(x)dx; g \in \mathbb{G} \}=\{ \widehat{g}(y); g \in \mathbb{G} \}; \ y \in \R^d.
\end{align}
\end{definition}

\begin{remark}
    These integrals can be extended to $\S'$-valued integrands by using the (weak) Gelfand-Pettis interpretation of the integral. See Remark \ref{rem1.7}. \\
    By Corollary 2.5.8 in \cite{HOUZ} we have the following:
    \begin{proposition}
        Suppose $Y(t,x)$ is a process such that
        \begin{align}
            \E[\int_{\R_{+}}\int_{\R^d} Y^2(t,x)dx dt ] < \infty.
        \end{align}
  Then $\int_{\R_{+}}\int_{\R^d} Y(t,x) \diamond W(t,x) dx dt$ exists as a Gelfand-Pettis integral, where $\diamond$ denotes the Wick product.\\  
  In particular, if in addition $Y(t,x)$ is deterministic , then
  $\int_{\R_{+}}\int_{\R^d} Y(t,x) W(t,x) dx dt$ exists as a Gelfand-Pettis integral.
  \end{proposition}
\end{remark}

\section{The fractional stochastic heat inclusion}
We are now ready to state and prove our main results. \\

The first part of the following result and its proof follow the approach in Theorem 2 in  \cite{MO1}.
\begin{theorem}
Let $G$ and $\Y$ be as in Section 1.1. \\
Let $\mathbb{G}=\mathbb{G}(Y)$ be the set of all bounded (and measurable) selectors of $G(t,x,Y(t,x))$, i.e. the set of all bounded processes $g$ such that $g(s,z) \in G(s,z,Y(s,z))\ \forall s,z$.

Then the following holds:
\begin{myenumerate}
    \item
    If $Y \in \Y$ is a  solution (in the sense of distribution) of the fractional stochastic heat inclusion
\begin{align}\label{heat}
 	 \frac{\partial^{\alpha}}{\partial t^{\alpha}}Y(t,x)&\in \l \Delta Y(t,x)+G(t,x,Y(t,x))+\sigma W(t,x);\nonumber\\
     &(t,x)\in (0,\infty)\times \mathbb{R}^{d},
 	 	\end{align}
then $Y$ is a fixed point of the \emph{integral inclusion}
\begin{align}
Y(t,x)&\in I_1(t,x) + I_2(t,x,Y(\cdot)) +I_3(t,x) \label{Y},
\end{align}
where
\small{
\begin{align}
I_1(t,x)=I&=(2\pi)^{-d} \int_{\mathbb{R}^d} e^{ixy}\sum_{k=0}^{\infty} E_{\alpha}(- \l t^{\alpha} |y|^2) dy\nonumber\\
=&(2\pi)^{-d} \int_{\mathbb{R}^d} e^{ixy}\sum_{k=0}^{\infty} \frac{(- \l t^{\alpha} |y|^2)^k}{\Gamma(\alpha k +1)}dy\label{I1}
\end{align}
\begin{align}
&I_2(t,x,Y(\cdot))
=\nonumber\\
&(2\pi)^{-d}\Big\{ \int_{0}^{t}(t-r)^{\alpha -1}
    \int_{\mathbb{R}^{d}}(\int_{\mathbb{R}^{d}}e^{i(x-\zeta)y}\sum_{k=0}^{\infty}\frac{(-\l (t-r)^{\alpha}y^2)^{k}}{\Gamma(\alpha (k+1))}dy) g(r,\zeta)dr d\zeta \Big\}_
     {g(\cdot) \in \G(Y)} \label{I2}
     \end{align}
     \begin{align}
    I_3(t,x)= (2\pi)^{-d} \int_{0}^{t}(t-r)^{\alpha -1}\int_{\mathbb{R}^{d}}\left(\int_{\mathbb{R}^{d}}e^{i(x-\zeta)y} E_{\alpha,\alpha}(-\l (t-r)^{\alpha}|y|^2) dy\right)
    \sigma B(dr,d\zeta)\nonumber\\
    =(2\pi)^{-d} \int_{0}^{t}(t-r)^{\alpha -1}\int_{\mathbb{R}^{d}}\left(\int_{\mathbb{R}^{d}}e^{i(x-\zeta)y}\sum_{k=0}^{\infty}\frac{(-\l (t-r)^{\alpha}|y|^2)^{k}}{\Gamma(\alpha k+\alpha))}dy\right)
\sigma B(dr,d\zeta). \label{I3}
\end{align}
}
\item
Conversely, if $Y \in \Y$ is a fixed point of the integral inclusion \eqref{Y}, then $Y$ is a solution of the fractional heat inclusion \eqref{heat}.
\end{myenumerate}
\end{theorem}
\dproof\\
In the following, we use the notation $|y|^2=y^2=\sum_{j=1}^d y_j^2; \quad y \in \mathbb{R}^d.$

\textbf{(i)} First assume that $Y(t,x)$ is a solution of \eqref{heat}.  We apply the Laplace transform $L$ 
to the $t$-variable both sides of  \eqref{heat} and obtain (see \eqref{L1}):
 \small
   \begin{align}
s^{\alpha}\widetilde{Y}(s,x)-s^{\alpha-1}Y(0,x)&\in \lambda \Delta \widetilde{Y}(s,x)+\Big \{ \widetilde{g}(s,x);  g(u,x)\in G(u,x,Y(u,x))\forall u) \Big\}\nonumber\\
    &+\sigma \widetilde{W}(s,x). 
	\end{align}
Applying the Fourier transform $F$, 
  	 we get, since $\widehat{Y}(0,y)=1$, 
  	 \small
     \begin{equation}
  	 	s^{\alpha}\widehat{\widetilde{Y}}(s,y)-s^{\alpha-1}\in -\l\sum_{j=1}^{d}y_{j}^{2}\widehat{\widetilde{Y}}(s,y)+\Big \{ \widehat{\widetilde{g}}(s,y);   g(u,z)\in G(u,z,Y(u,z))\forall u,z)\Big\}+\sigma\widehat{\widetilde W}(s,y).
  	 \end{equation}  
   Hence
   \begin{align}
       \widehat{\widetilde{Y}}(s,y)\in \frac{s^{\alpha -1}}{s^{\alpha} + \l y^2}
       +\Big \{ \frac{\widehat{\widetilde{g}}(s,y)}{s^{\alpha} +\l y^2};   g\in\G(Y) \Big\}+\frac{\sigma \widehat{\widetilde{W}}(s,y)}{s^{\alpha}+\l y^2}.
   \end{align}
  Since the Laplace transform and the Fourier transform commute, this can be written
   \begin{align}
       \widetilde{\widehat{Y}}(s,y)\in \frac{s^{\alpha -1}}{s^{\alpha} + \l y^2}
       +\Big \{ \frac{\widetilde{\widehat{g}}(s,y)}{s^{\alpha} + \l y^2};   g\in\G(Y) \Big\} +\frac{\sigma \widetilde{\widehat{W}}(s,y)}{s^{\alpha}+\l y^2}.
   \end{align}
 Applying the inverse Laplace operator $L^{-1}$  we get 
\begin{align}
	L^{-1}\left(\frac{1}{s^{\alpha}+\l y^2} \right) (t)&=t^{\alpha-1}E_{\alpha,\alpha}(-\l t^{\alpha}y^2)
	=\sum_{k=0}^{\infty}\frac{t^{\alpha-1}(-\l t^{\alpha}y^2)^{k}}{\Gamma(\alpha k+\alpha)}\nonumber\\
	&=\sum_{k=0}^{\infty}\frac{(-\lambda y^2)^{k}t^{\alpha(k+1)-1}}{\Gamma(\alpha(k+1))}
	=\sum_{k=0}^{\infty}\frac{(-\l t^{\alpha}y^2)^{k}t^{\alpha -1}}{\Gamma(\alpha (k+1))}	
	=: \Lambda(t,y).\label{Lambda}
\end{align}
In other words,
\begin{equation}
	\frac{1}{s^{\alpha}+\l y^2}= L \Lambda(s,y),
\end{equation}
Hence
\small
 \begin{align} 
     & \widehat{Y}(t,y)\nonumber\\
     &\in E_{\alpha,1}(-\l |y|^2 t^{\alpha})+L^{-1}\Big \{ \frac{\widetilde{\widehat{g}}(s,y)}{s^{\alpha} + \l y^2};   g\in\G(Y)  \Big\}  + L^{-1}\Big(\frac{\sigma \widetilde{\widehat{W}}(s,y)}{s^{\alpha}+\l y^2}\Big)(t,y),\nonumber\\
       &= E_{\alpha,1}(-\l |y|^2 t^{\alpha})+\Big \{ L^{-1}\Big(\frac{\widetilde{\widehat{g}}(s,y)}{s^{\alpha} + \l y^2}\Big);   g\in\G(Y)\Big\}  + L^{-1}\Big(\frac{\sigma \widetilde{\widehat{W}}(s,y)}{s^{\alpha}+\l y^2}\Big)(t,y),\label{8}\nonumber\\
 &= E_{\alpha,1}\left( -\l t^{\alpha}y^2\right)+\Big \{ \int_0^t \Lambda(t-r,y)\widehat{g}(r,y)dr;   g\in\G(Y) \Big\}\nonumber\\
 &+\int_{0}^{t} \Lambda(t-r,y)\sigma \widehat{W}(r,y)dr.
\end{align}
Taking inverse Fourier transform we end up with
\small{
\begin{align}\label{2.18}
 &Y(t,x)\nonumber\\
 &\in F^{-1}\left( E_{\alpha,1}\left(-\l t^{\alpha}|y|^2 
 \right)\right)(x)+F^{-1}\Big \{ \int_0^t \Lambda(t-r,y)\widehat{g}(r,y)dr;   g\in\G(Y) \Big\}(x)\nonumber\\
 &+ F^{-1}\left(\int_{0}^{t}\Lambda(t-r,y)\sigma \widehat{W}(r,y)dr\right)(x).
\end{align}

This can be written
 \small
\begin{align*}
   & Y(t,x)\nonumber\\
    &\in (2\pi)^{-d} \int_{\mathbb{R}^d} e^{ixy}\sum_{k=0}^{\infty} \frac{(- \l t^{\alpha} y^2)^k}{\Gamma(\alpha k +1)}dy \nonumber\\
    &+  (2\pi)^{-d}\Big\{ \int_{0}^{t}(t-r)^{\alpha -1}
    \int_{\mathbb{R}^{d}}\left(\int_{\mathbb{R}^{d}}e^{i(x-z)y}\sum_{k=0}^{\infty}\frac{(-\l (t-r)^{\alpha}y^2)^{k}}{\Gamma(\alpha (k+1))}dy\right) g(r,z)dr dz ;  g\in\G(Y) \Big\}\\
    &+ (2\pi)^{-d} \int_{0}^{t}(t-r)^{\alpha -1}
    \int_{\mathbb{R}^{d}}\left(\int_{\mathbb{R}^{d}}e^{i(x-z)y}\sum_{k=0}^{\infty}\frac{(-\l (t-r)^{\alpha}y^2)^{k}}{\Gamma(\alpha (k+1))}dy\right)\sigma B(dr,dz) \\
    &=I_1(t,x) + I_2 (t,x,Y(\cdot))+ I_3(t,x)
\end{align*}

This proves that the solution $Y \in \Y$ (if it exists) is a fixed point of the inclusion \eqref{Y},
as claimed.

\textbf{(ii)} This follows by reversing the argument in (i). We skip the details.

\begin{remark}
    Note that the fixed point inclusion \eqref{Y} is of the form $Y(t,x) \in H(Y(\cdot))$ for some set valued function $H$, where the right hand side depends on the whole path $\{Y(s,a); s \leq t, a \in \R^d \}$ up to $t$, and not on just the terminal value $Y(t,x)$. Therefore this fixed point inclusion is not of the usual type and it is not covered by the fixed point inclusion theorems known to us. Because of the dependence of the whole path  on its right hand side,  \eqref{Y} could be called a \emph{Volterra inclusion}, in analogy of the situation for classical Volterra integral equations.
\vskip 0.2cm
    Proving an existence and uniqueness theorem for Volterra inclusions would be a challenging new project. We will not discuss it further 
    here.
 \end{remark}
 \bigskip

In the special case when the drift set $G$ does not depend on $Y$, we no longer have fixed point problem:
\begin{corollary}
Assume that $G(t,x,Y(t,x))=G(t,x)$ does not depend on $Y(t,x)$. Let
$\mathbb{G}$ be the set of all bounded (and measurable) selectors of $G(t,x)$, i.e. the set of all bounded processes $g$ such that $g(t,x) \in G(t,x)\ \forall t,x$.\\
Then the solution of the stochastic differential inclusion \eqref{heat} is given directly by the following inclusion:
\begin{align}
Y(t,x)&\in I_1(t,x) + I_2(t,x) +I_3(t,x) \label{Y2},
\end{align}
where
\begin{align}
I_1(t,x)&=(2\pi)^{-d} \int_{\mathbb{R}^d} e^{ixy}\sum_{k=0}^{\infty} E_{\alpha}(- \l t^{\alpha} |y|^2) dy\nonumber\\
=&(2\pi)^{-d} \int_{\mathbb{R}^d} e^{ixy}\sum_{k=0}^{\infty} \frac{(- \l t^{\alpha} |y|^2)^k}{\Gamma(\alpha k +1)}dy\label{I12},
\end{align}
\begin{align}
&I_2(t,x)
=\nonumber\\
&(2\pi)^{-d}\Big\{ \int_{0}^{t}(t-r)^{\alpha -1}
    \int_{\mathbb{R}^{d}}(\int_{\mathbb{R}^{d}}e^{i(x-\zeta)y}\sum_{k=0}^{\infty}\frac{(-\l (t-r)^{\alpha}y^2)^{k}}{\Gamma(\alpha (k+1))}dy) g(r,\zeta)dr d\zeta \Big\}_
     {g(\cdot) \in \G}, \label{I22}
     \end{align}
     \begin{align}
    I_3(t,x)= (2\pi)^{-d} \int_{0}^{t}(t-r)^{\alpha -1}\int_{\mathbb{R}^{d}}\left(\int_{\mathbb{R}^{d}}e^{i(x-\zeta)y} E_{\alpha,\alpha}(-\l (t-r)^{\alpha}|y|^2) dy\right)
    \sigma B(dr,d\zeta)\nonumber\\
    =(2\pi)^{-d} \int_{0}^{t}(t-r)^{\alpha -1}\int_{\mathbb{R}^{d}}\left(\int_{\mathbb{R}^{d}}e^{i(x-\zeta)y}\sum_{k=0}^{\infty}\frac{(-\l (t-r)^{\alpha}|y|^2)^{k}}{\Gamma(\alpha k+\alpha))}dy\right)
\sigma B(dr,d\zeta). \label{I32}
\end{align}
\end{corollary}

\section{When is the solution mild?}

In the previous section, we proved that there is at least one distribution-valued solution $Y \in \Y.$ As mentioned in the Introduction, we define a solution $Y$ as \emph{mild} if $Y(t,x) \in L^2(\P)$ for all $t,x$.
 We now ask:
 
\emph{Are any of these solutions mild?}

This question has been discussed in the fractional heat \emph{equation} case in \cite{MO1}. 
In the following, we recall Theorem 3 in \cite{MO1}:
\begin{theorem}\cite{MO1}
 Let \( Y(t, x) \) be the solution of the \( \alpha \)-fractional stochastic heat equation given by:
 \begin{equation}\label{beya}
     \frac{\partial^\alpha}{\partial t^\alpha} Y(t, x) = \lambda \Delta Y(t, x) + \sigma W(t, x), \quad (t, x) \in (0, \infty) \times \mathbb{R}^d.
 \end{equation}
Then the following holds:  
\begin{myenumerate}
\item If \( \alpha = 1 \), then \( Y(t, x) \) is mild if and only if \( d = 1 \).  
 \item If \( \alpha > 1 \), then \( Y(t, x) \) is mild if \( d = 1 \) or \( d = 2 \).  
 \item  If \( \alpha < 1 \), then \( Y(t, x) \) is not mild for any \( d \).
 \end{myenumerate}
\end{theorem}
Return now to our case. By equation \eqref{Y}, we have \begin{equation}\label{Y3}
    Y(t,x) \in I_1(t,x)+I_3(t,x)+I_2(t,x,Y(.)).
\end{equation}
We want to study when the right-hand side of equation \eqref{Y3} belongs to $L^2(\P)$. 

Note that $I_1(t,x)+I_3(t,x)$ of equation \eqref{Y3} corresponds to the same term as in Theorem 3 in \cite{MO1}. 
Using Theorem 3  in \cite{MO1}, we get $\E[(I_1(t,x) +I_3(t,x))^2]=\E[C^2(t,x)] < \infty$ in the following cases:
\begin{myenumerate}
    \item 
    If $\alpha=1$ and $d=1$
    \item 
    If $\alpha > 1$ and $d=1$ or $d=2$
\end{myenumerate}
Moreover, it is proved that if $\alpha <1$ then $\E[C^2(t,x)] = \infty$ for all $d$.\\
It remains to consider 
\begin{equation}
    I_2(t,x)=:\int_0^t \int_{\R^d} h_0(t,x;s,z)g(s,z)ds dz, \quad\text{  for all  } g \in G(Y)
\end{equation}
in the case when $\alpha >1$.
Here $h_0$ is given by
 \begin{align}
     h_0(t,x;s,z)=(t-s)^{\alpha -1}\int_{\R^d} e^{i(x-z)y}\sum_{k=0}^{\infty}\frac{(-\l (t-s)^{\alpha}y^2)^{k}}{\Gamma(\alpha (k+1))}dy.\label{g}
    \end{align}
    Note that
    \begin{align}
        I_3(t,x)=(2\pi)^{-d} \int_0^t \int_{\R^d} h_0(t,x; s,z)\sigma B(ds, dz). \label{g2}
    \end{align}

By the same argument as in Theorem 3b) in \cite{MO1} we get that for all $g \in G(Y)$ we have, for $\alpha >1$ and $d=1$ or $d=2$,
\begin{align}
    &\E\Big[\Big(\int_0^t\int_{\R^d} h_0(t,x;s,\zeta) g(s,\zeta) ds d\zeta \Big)^2 \Big]\nonumber\\
    &\leq \E[\int_{S(g)}g^2(s,\zeta)dsd\zeta] \  \E[\int_0^t \int_{\R^d} h_0^2(t,x;s,\zeta) ds d\zeta] \nonumber\\
    & \leq const. \  \E[\int_0^t \int_{\R^d} h_0^2(t,x;s,\zeta) ds d\zeta] < \infty.
\end{align}

Moreover, if $\alpha >1$ and $d=1$ or $d=2$ we know that $\E[C^2(t,x)] < \infty$ and we conclude that $\E[Y^2(t,x)] = \infty$.\\
We summarize the above as follows:
\begin{theorem} \label{th8.1}
\textbf{a)}   Suppose one of the following holds:\begin{myenumerate}
    \item 
    $\alpha=1$ and $d=1$
    \item 
    $\alpha > 1$ and $d=1$ or $d=2.$
\end{myenumerate}
Then there exists a mild solution of the fractional heat inclusion \eqref{heat}.\\
\textbf{b)} If $\alpha <1$ there is no mild solution of \eqref{heat}.
\end{theorem}

\section{Example: Nitrate transport in sandy porous aquifers}\label{new}

We end this paper with an example, illustrating how our theoretical result  which combines fractional derivatives, set-valued maps and white noise can be used to describe the spreading of nitrate pollution through sandy groundwater systems.

Sandy aquifers are loose or dense sand that is water-permeable and dissolved
substances to pass through tiny pores between the grains. Such systems occur in
vast portions of the world and are mostly vulnerable to poisoning due to agricultural
activities, such as overuse of fertilisers \cite{fetter2001}.

We focus on the transport of nitrate, our model is based on the following fractional stochastic inclusion:

$\frac{\partial^\alpha}{\partial t^\alpha} Y(t, x) \in \lambda \Delta Y(t, x) + G(Y(t, x)) + \sigma W(t, x); \quad (t, x) \in (0, \infty) \times \mathbb{R}^d,$\\
with initial and boundary conditions:
$$
Y(0, x) = \delta(x), \quad \lim_{|x| \to \infty} Y(t, x) =0
$$
Let us modelise term by term:

The fractional time derivative $\frac{\partial^\alpha}{\partial t^\alpha} Y(t, x)  $ captures memory effects and anomalous diffusion.

When $\alpha<1 $, it models slower spreading and trapping of nitrate.

When $\alpha>1$, it reflects faster-than-normal movement, such as during sudden water infiltration.

The diffusion term  describes how nitrate naturally spreads out due to differences in  concentration.

The linear set-valued map $G(Y)=[k_1 Y,k_2 Y]$, where $0 \leq k_1 < k_2$ are given constants, represents uncertainty in nitrate absorption, depending on local soil properties. Sandy soils might absorb very little (low $k_1$), while areas with organic matter might absorb more (higher $k_2$) \cite{fetter2001}.

The white noise term $W$  models random microscopic variations in permeability and flow conditions.

About the initial condition:

The first means that there   is a progressive release of nitrate from a source at a location, such as
a flowing stream, continuous seepage, or flowing aquifer. This tells us about
how pollution builds up from a diffuse source.

In the real world, this might equate to an initial concentration  of 50 to $200 mg/L$, \cite{wakida2005},  much higher
than the World Health Organization suggested limit of $50 mg/L$ for potable drinking water \cite{who2017}.

This model is especially useful for simulating how nitrate travels through sandy aquifers, where its high permeability and low absorption combine to make it a very real threat. The fractional-time derivative allows us to simulate the fact that nitrate does not travel uniformly at a steady rate; it may become stuck, then zip ahead suddenly, depending on soil texture and water movement.
  
  Brownian noise terms and set-valued uncertainty allows us   to simulate the effects of unknown parameters, natural variability, and unexpected events, all of which happen in real aquifer systems.\\
Now, question is how our models can assist and predict the risk of a high concentration of nitrate.

Our solution $Y(t,x)$ is the the logarithm of the concentration of nitrate at time t
 
and position x, it captures the way the nitrate propagates through fractures or porous media.
It also shows the way that the system retains past states because of the fractional derivative. The high concentration of nitrate in the groundwater has extreme perils, mostly
to human health and to the environment. For example, we bring to light the Blue baby
syndrome, which infects babies by reducing their blood's oxygen-carrying ability.
Aside from that, there is also growing evidence of the presence of a link between long-term exposure to nitrate and certain cancers, especially stomach cancer.

The solution $Y(t,x)$ tells us where and when the concentration  of nitrates would be above the safe level. It also tells us the regions of potential danger zones ahead, regions where the nitrate concentration  could become excessive. It helps decision-makers and communities prepare by telling them where and when protective actions like filter installation would be needed, and it reveals which places are most vulnerable and why.
\subsection{Numerical Example}

We present in this section a numerical simulation of the following fractional stochastic differential inclusion:

\[
\frac{\partial^\alpha}{\partial t^\alpha} Y(t, x) \in \lambda \Delta Y(t, x) + G(Y(t, x)) + \sigma W(t, x)  , \quad (t, x) \in (0, \infty) \times \mathbb{R}.
\]

The simulation uses the following parameters:
\begin{itemize}
    \item Fractional order: \( \alpha = 0.8 \).
    \item Diffusion coefficient: \( \lambda = 0.1 \).
    \item Noise intensities: \( \sigma = 0.5 \).
    \item  \( G(Y) = [0.1Y, 0.5Y] \), representing a range of absorption intensities.
    \item Domain: spatial range \( x \in [-50, 50] \) meters, over a time horizon \( t \in [0, 30] \) days.
\end{itemize}
\begin{figure}[h!]
    \centering
    \includegraphics[width=0.95\linewidth]{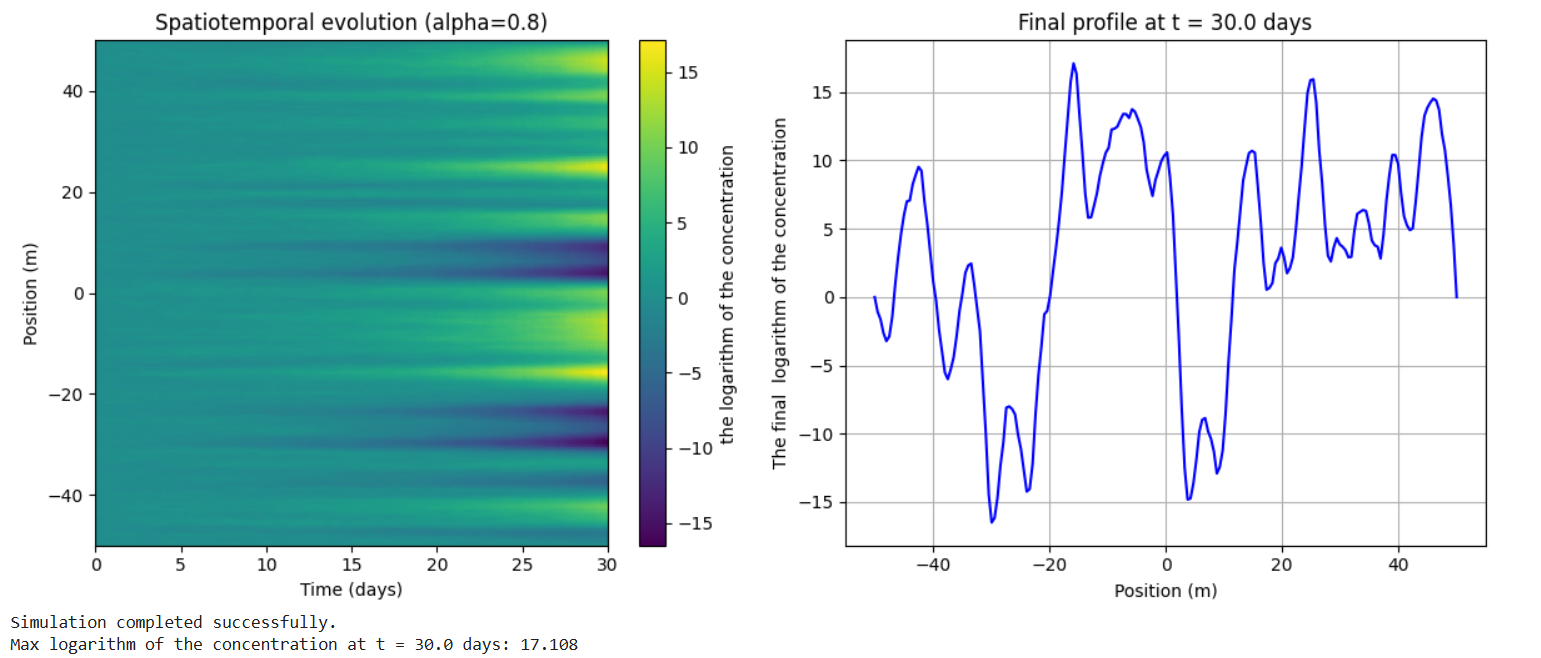}

        \caption{Spatiotemporal evolution of the  logarithm of the concentration \( Y(t,x) \).}
    \label{fig:spatiotemporal}
\end{figure}

\textbf{Interpretation of the results of Figure 1: }
\paragraph{Left Plot:}

\textbf{What It Shows:}
\begin{itemize}
    \item X-axis: Time.  
    \item Y-axis: Spatial position.
    \item Color:  log concentration at each position over time
    \item Dark blue: Low log concentration 
    \item Bright yellow: High log concentration 
\end{itemize}
We interpret the solution $Y(t,x)$  as the logarithm of the concentration of the nitrate, then negative values for  $Y(t,x)$ simply correspond to concentrations between 0 and 1.
  The threshold becomes $Y(t,x)=0$, which corresponds to a concentration of 1 unit. We would then define the "danger zone" as the region where $Y(t,x)>0$ . We can calibrate this "1 unit" to represent the World Health Organization safety limit of 50 mg/L, meaning we set our model such that  log(50) is the new zero point, or simply state that the critical threshold is 0 on our log-scale.
   
\begin{remark}
    It is natural to use a logarithmic scale when data spans a very wide range of values, when showing exponential growth or trends, or when comparing percentage changes. Using a log scale allows us to see both small and large values clearly on the same chart, which would otherwise be compressed or unreadable on a linear scale. It is common in scientific fields like medicine and archaeology, finance, and when visualizing data with a high degree of skewness. 
\end{remark}
 
 \subsection{Tools for the numerical method}
We consider the stochastic fractional diffusion-absorption inclusion
\begin{equation}
\partial_t^\alpha Y(t,x)
= \lambda \Delta Y(t,x)
  - k(t) Y(t,x)
  + \sigma \dot{W}(t,x),
\qquad (t,x) \in [0,T] \times \Omega,
\label{eq:model}
\end{equation}
where $0<\alpha<1$ is the Caputo fractional order, 
$\lambda>0$ is the diffusion coefficient, 
$k(t)\in[k_1,k_2]$ is a randomly varying absorption rate,
and $\dot{W}$ denotes Gaussian space--time white noise.
\subsubsection{Spatial discretization}
The spatial domain $\Omega=[-L/2,L/2]$ is discretized using a uniform grid
$x_i=-L/2+i\Delta x$, $i=0,\dots,N_x$, with $\Delta x = L/N_x$.  
The Laplacian is approximated by the standard centered second-order finite difference:
\[
(\Delta Y)^n_i
\approx 
\frac{Y^n_{i+1}-2Y^n_i+Y^n_{i-1}}{\Delta x^2}.
\]
Dirichlet boundary conditions are imposed by setting
$Y^n_0 = Y^n_{N_x}=0$ for all $n$.
\subsubsection{Time discretization}
The time interval $[0,T]$ is discretized uniformly with 
$t_n = n\Delta t$, $n=0,\dots,N_t$, where $\Delta t = T/N_t$.
The update is explicit:
\[
Y^{n}_i
=
Y^{n-1}_i
+ \Delta t
\Big(
\lambda (\Delta Y)^{n-1}_i
- k^n Y^{n-1}_i
- M^{n}_i
\Big)
+ \sigma \sqrt{\Delta t}\,\eta^n_i,
\]
where $\eta^n_i\sim\mathcal{N}(0,1)$ are i.i.d.\ Gaussian variables.
\subsubsection{Fractional Caputo derivative approximation}
The Caputo derivative is approximated using a truncated convolution sum:
\[
\partial_t^\alpha Y(t_n,x_i)
\approx
\frac{1}{\Gamma(2-\alpha)\Delta t^\alpha}
\sum_{j=1}^{m^*}
\big(Y^{n-j}_i - Y^{n-j-1}_i\big)
\left( j^{1-\alpha} - (j-1)^{1-\alpha} \right),
\]
where $m^*$ is the chosen memory length (here $m^*=15$).  The quantity
\[
M^n_i
= 
\frac{1}{\Gamma(2-\alpha)\Delta t^\alpha}
\sum_{j=1}^{m^*}
\big(Y^{n-j}_i - Y^{n-j-1}_i\big)
\left( j^{1-\alpha} - (j-1)^{1-\alpha} \right)
\]
acts as a history-dependent drift term.
\subsubsection{Random nonlinear absorption}
At each time step, the absorption coefficient is sampled uniformly:
\[
k^n \sim \mathcal{U}[k_1,k_2],
\qquad
g(Y^{n-1}_i) = k^n Y^{n-1}_i.
\]
This produces a nonlinear, randomly switching absorption term,  
thereby turning the equation into a stochastic differential inclusion.
\subsubsection{Noise term}
Space--time white noise is approximated explicitly by:
\[
\sigma \sqrt{\Delta t}\,\eta_i^n,
\qquad
\eta_i^n \sim \mathcal{N}(0,1).
\]
\subsubsection{Summary}
The global numerical scheme is an explicit Euler-type method with:
\begin{itemize}
    \item finite-difference spatial discretization,
\item Grunwald-type memory approximation for the Caputo derivative,
    \item random multivalued absorption $k(t) \in [k_1,k_2]$,
    \item additive Gaussian noise.
\end{itemize}
The method produces the spatiotemporal solution matrix
\[
Y^n_i \approx Y(t_n,x_i),
\qquad
n=0,\dots,N_t,\; i=0,\dots,N_x,
\]
from which all plots are generated.
}
\begin{remark}
   While Theorem \ref{th8.1} establishes the non-existence of mild solutions when $\alpha < 1$, our numerical results  demonstrate stable approximations. This   regularity come from the following aspects:

\begin{enumerate} 
    \item \emph{Domain truncation}: The finite spatial domain $x \in [-L,L]$ and time horizon $T < \infty$ avoid the non-integrability of $E_{\alpha,\alpha}(-t^\alpha)$ at infinity.
    
    \item \emph{Kernel  }: The convolution with the Mittag-Leffler function replaces the singular kernel with a bounded approximation:
    \begin{equation*}
        \widetilde{E}_{\alpha,\alpha}(t) := 
        \begin{cases}
            E_{\alpha,\alpha}(-t^\alpha) & \text{for } t \geq \Delta t \\
            0 & \text{otherwise}
        \end{cases}
    \end{equation*}
    
    \item  

Mathematically,  the numerical solution approximates a \emph{weak solution} in $\mathcal{S}'(\mathbb{R}^d)$: For any test function $\phi \in C_c^\infty(\mathbb{R}^d)$, the discrete scheme satisfies
\begin{equation}
    \lim_{\substack{\Delta t \to 0 \\ \Delta x \to 0}} \mathbb{E}\left[\left|\langle Y^{\Delta t,\Delta x}_n, \phi\rangle - \langle Y, \phi\rangle\right|^2\right] = 0
\end{equation}
where $\langle \cdot, \cdot \rangle$ denotes the duality pairing. This aligns with the regularizing effects observed in stochastic PDEs with multiplicative noise \cite[Theorem 3.5.1]{oksendal1996} and fractional dynamics \cite[Section 4]{MS}.
  \end{enumerate}
\end{remark}


\begin{thebibliography}{99}
   
   \bibitem{A} Abel, N. H. (1823): Oppl\o sning av et par oppgaver ved hjelp av bestemte integraler (in Norwegian).  Magazin for naturvidenskaberne 55--68.

\bibitem{AMSW} Adler, R.J., Monrad, D., Scissors, R.H. \& Wilson, R. (1983): Representations, decompositions and sample
function continuity of random fields with independent increments. Stoch. Process. Appl. 15(1),3-30.

  \bibitem{Au} Aumann, R. J. (1965): Integrals of set-valued functions. J. Math.Anal. 12, 1-12.
  
\bibitem{CMV} Capelas de Oliveira, E.; Mainardi \& F.Vaz, J. (2011): Models based on Mittag Leffler functions for anomalous relaxation in dielectrics. ariXir: 1106.1761  v2[cond-mat.Stat-mech] 13 Feb. 2014.

\bibitem{DW} Dalang, R. C. \& Walsh, J. B. (1992): The sharp Markov property of L\'evy sheets. The Annals of Probability 20 (2), 591-626.

\bibitem{DZ}Da Prato, G. \& Zabczyk, J. (2014): Stochastic Equations in Infinite Dimensions.
Second edition. Encyclopedia of Mathematics and its Applications, 152. Cambridge
University Press, Cambridge.
  \bibitem{DOP} Di Nunno, G., \O ksendal, B. \& Proske, F. (2009): Malliavin Calculus for L\' evy Processes with Applications to Finance. Springer.
  \bibitem{fetter2001}
Fetter, C.~W. (2001). \textit{Applied Hydrogeology} (4th ed.). Prentice Hall.
   \bibitem{Hu} Hu, Y. (2019): Some recent progress on stochastic heat equations. Acta Mathematica Scientia 39B(3); 874-914.
   
   \bibitem{H} Holm, S. (2019): Waves with Power-Law Attenuation. Springer.


 \bibitem{HOUZ}  Holden, H., \O ksendal, B.,  Ub\o e, J.  and Zhang, T.: Stochastic Partial Differential Equations.
Birkh\"auser Boston Inc., Boston, MA, 1996.

 \bibitem{I} Ibe, O. C. (2013): Markov Processes for Stochastic Modelling. 2$^{nd}$ edition. Elsevier. 

 \bibitem{IR} Iafrate, F. \& Ricciutu, C.(2024): Some families of random fileds related to multiparameter L\' evy processes. J. Theoretical Probability 37: 3055-3088. https://doi.org/10.1007/s10959-024-01351-3.
 
 \bibitem{GV}I.M. Gel'fand and N.Ya. Vilenkin. Generalized Functions. Vol. 4. Academic Press [Harcourt Brace Jovanovich Publishers], New York, 1964 [1977].


   \bibitem{Kilbas} Kilbas, A.A., Srivastava,H.M. \& Trujillo, J. J. (2006): 
	 {\emph{ Theory and Applications of Fractional Differential Equations }}, Elsevier  Science B.V, Amsterdam. 
     
\bibitem{KKS} Kochubel, A. N., Kondratiev, Y. \& da Silva, J. L. (2021): On fractional heat equation. Fractional Calculus \& Applied Analysis 24 (1), 73-87.
	     
   \bibitem{MS} Meerschaert, Mark M., Sikorskii, Alla (2019): Stochastic Models for Fractional Calculus, $2^{nd}$ edition. De Greuter.
   
   \bibitem{ML} Mittag-Leffler, M.G. (1903): Sur la nouvelle fonction E(x). \emph{C. R. Acad. Sci. Paris} {\bf 137} 554--558.
   
  \bibitem{MO1} Moulay Hachemi, R.Y.  \& \O ksendal, B. (2023): The fractional stochastic heat equation
driven by time-space white noise. Fractional Calculus \& Applied Analysis 26, 513-532. Springer.
https://doi.org/10.1007/s13540-023-00134-7.

\bibitem{NualartSchoutens} Nualart, D. \& Schoutens, W. (2000): Chaotic and predictable representations for L\'e vy processes.
Stochastic Process. Appl., 90(1):109-122, 2000.

\bibitem{oksendal1996}
\O ksendal, B. (2013):  Stochastic Differential Equations, 6th edition, Springer, Heidelberg.

\bibitem{OP} \O ksendal, B. and Proske, F. (2004): White noise of Poisson random measures. Potential Anal., 21(4):375-403.

\bibitem{PK} Patel, A. \& Kosko, B. (2007): "Levy Noise Benefits in Neural Signal Detection," 2007 IEEE International Conference on Acoustics, Speech and Signal Processing - ICASSP '07, Honolulu, HI, USA, 2007, pp. III-1413-III-1416, doi: 10.1109/ICASSP.2007.367111. 

\bibitem{P} Pollard, H. (1948): The completely monotone character of the Mittag-Leffler function $E_{\alpha}(-x)$. Bull. Amer. Math-Soc.54, 1115-1116. 



\bibitem{Samko} Samko, S.G., Kilbas, A.A., Marichev, O.I. (1993): Fractional Integrals and Derivatives. Theory and Applications. Gordon and Breach Science Publishers, New York.

\bibitem{S} Schneider, W.B. (1996): Completely monotone generalized Mittag-Leffler functions. Expositiones Mathematicae 14, 3-16.

\bibitem{wakida2005}
Wakida, F.~T. and Lerner, D.~N. (2005). Non-agricultural sources of groundwater nitrate: a review and case study. \textit{Water Research}, 39(1):3--16.
\bibitem{who2017} 
World Health Organization (2017). \textit{Guidelines for Drinking-water Quality: Fourth Edition Incorporating the First Addendum}. WHO Press, Geneva.
\bibitem{Aubin2003} Aubin, J.-P., Frankowska, H. (2009). \textit{Set-Valued Analysis}. Springer Science and Business Media.

\bibitem{glicksberg1952} Glicksberg, I. L. (1952). A Further Generalization of the Kakutani Fixed Point Theorem, with Application to Nash Equilibrium Points. \textit{Proceedings of the American Mathematical Society}, 3(1), 170.174.
 
\bibitem{Rudin1991} Rudin, W. (1991). \textit{Functional Analysis} (2nd ed.). McGraw-Hill.
\bibitem{aubin2009} Aubin, J.-P., Frankowska, H. (2009). \emph{Set-Valued Analysis.}
\bibitem{zeidler1986} Zeidler, E. (1986). \emph{Nonlinear Functional Analysis and its Applications I.}
\bibitem{podlubny1998} Podlubny, I. (1998). \emph{Fractional Differential Equations.}
\bibitem{da2014} Da Prato, G., Zabczyk, J. (2014). \emph{Stochastic Equations in Infinite Dimensions.}
\bibitem{fan1952} K. Fan, Fixed-point and minimax theorems in locally convex topological linear spaces, Proceedings of the National Academy of Sciences, 38(2):121-126, 1952.
\end{thebibliography}
\end{document}